\documentclass[11pt]{article}

\usepackage{latexsym}
\usepackage{amsmath}
\usepackage{amssymb}
\usepackage{amsthm}
\usepackage{hyperref}
\usepackage{cite}
\usepackage{graphicx}
\usepackage{color}
\usepackage{enumitem}
\usepackage{braket}
\usepackage{bm}
\usepackage[left=1in,right=1in,top=1in,bottom=1in]{geometry}
\usepackage{xspace}
\setlist[description]{font=\normalfont\itshape\textbullet\space}
\usepackage[all]{xy}

\renewcommand{\paragraph}[1]{\vspace{6pt} \noindent \textbf{#1}\xspace}

\theoremstyle{plain}
\newtheorem{theorem}{Theorem}[section]

\newtheorem{lemma}[theorem]{Lemma}

  {\end{tabular}\par\medskip}

\newtheorem{postulate}{Postulate}

\theoremstyle{definition}

\newtheorem{definition}[theorem]{Definition}
\newtheorem{example}[theorem]{Example}

\newcommand{\GL}{\mathrm{GL}}
\newcommand{\F}{\mathbb{F}}
\newcommand{\Z}{\mathbb{Z}}

\newcommand{\GalN}[2]{\mathrm{GalN}^{#1}_{#2}}
\newcommand{\GN}{\mathrm{G}}
\newcommand{\NDG}{\mathrm{NDG}}

\newcommand{\NDS}{\mathrm{NDS}}
\newcommand{\RGN}{\mathrm{RG}}
\newcommand{\CGN}{\mathrm{CG}}
\newcommand{\DIS}{\mathrm{DIS}}
\newcommand{\RSN}{\mathrm{RS}}

\newcommand{\N}{\mathbb{N}}

\newcommand{\qbinom}[3]{\genfrac{[}{]}{0pt}{}{#1}{#2}_{#3}}

\newcommand{\rad}{\mathrm{rad}}

\newcommand{\M}{\mathrm{M}}

\newcommand{\spa}[1]{\mathcal{#1}}

\newcommand{\cA}{\spa{A}}
\newcommand{\cB}{\spa{B}}

\newcommand{\aut}{\mathrm{Aut}}

\newcommand{\too}%
{\xrightarrow{\text{\raisebox{-3pt}{$\sim$}}\,}}

\usepackage{arydshln}
\hypersetup{
    pdftitle={Enumerating alternating matrix spaces},
    pdfauthor={Youming Qiao},
    bookmarksnumbered=true,     
    bookmarksopen=true,         
    bookmarksopenlevel=1,       
    colorlinks=true,            
    pdfstartview=Fit,           
    pdfpagemode=UseOutlines,    % this is the option you were lookin for
    pdfpagelayout=OneColumn,
    pdfstartview=FitH
}
\title{
Enumerating alternating matrix spaces over finite fields with explicit coordinates 
}

\author{
Youming Qiao
\footnote{Centre for Quantum Software and Information, University of 
Technology Sydney. \tt{jimmyqiao86@gmail.com}} 
}
\date{\today}

\begin{document}

%\pagenumbering{gobble}  % no page number on first page

\maketitle

\begin{abstract} 
We initiate the study of enumerating linear subspaces of 
alternating matrices over finite fields with explicit coordinates. 
We postulate
that this study can be viewed as a linear algebraic analogue of the classical 
topic of enumerating 
labelled graphs. 
To support this viewpoint, we present 
$q$-analogues of Gilbert's formula for enumerating connected graphs (\emph{Can. J. 
Math.}, 1956), and Read's formula for enumerating 
$c$-colored graphs (\emph{Can. J. Math.}, 
1960). We also develop an analogue of Riddell's formula relating the exponential 
generating function of 
graphs with that of connected graphs (Riddell's PhD thesis, 1951), building on 
Eulerian generating functions developed by Srinivasan (\emph{Discrete Math.}, 
2006).
\vskip 1em
\noindent{\it Keywords: labelled graph enumeration; alternating matrix spaces; 
Eulerian generating functions; q-calculus; alternating bilinear maps. }
%On the one hand, this study can be viewed as a linear algebraic analogue of 
%enumerating labelled graphs. On the other hand, 
%We put forward 
%the postulate 
%that this study can be 
%viewed as a linear algebraic analogue of the classical topic of enumerating 
%labelled graphs in graph theory. 
\end{abstract}

\section{Introduction}

\subsection{From enumerating graphs to enumerating alternating matrix spaces}

An $n\times n$ matrix $A$ over a field $\F$ is \emph{alternating}, if 
for any $v\in \F^n$, $v^tAv=0$. 
%We shall focus on alternating matrix spaces over 
%finite fields. 
Let $\Lambda(n, \F)$ be the linear space of $n\times n$ alternating 
matrices over $\F$. A subspace of $\Lambda(n, \F)$ 
is referred to as an 
\emph{alternating matrix space}. When $\F$ is the finite field with $q$ elements 
$\F_q$, we may also write $\Lambda(n, 
\F_q)$ as $\Lambda(n, q)$.

We study enumerating alternating matrix spaces over finite fields with explicit 
coordinates, and propose the following postulate.
\begin{postulate}\label{main}
Enumerating alternating matrix spaces over finite fields with explicit coordinates 
can be viewed and studied as a linear algebraic analogue of enumerating 
labelled graphs. 
\end{postulate}

To get a first hint of the analogy, recall that the number of labelled graphs with 
$n$ vertices and $m$ edges is $\binom{\binom{n}{2}}{m}$. Then note that the number 
of $m$-dimensional alternating matrix spaces in $\Lambda(n, q)$ is 
$\qbinom{\binom{n}{2}}{m}{q}$, 
where the $\qbinom{\ }{\ }{q}$ denotes the Gaussian binomial coefficient. 

One may wonder what this easy example could lead to. For this, we
reflect on 
the research into enumerating graphs, a highly productive research line, as marked 
by e.g. the classical monograph of Harary and Palmer \cite{HP73} and the recent 
wonderful survey of Wormald \cite{Wor18}. One key reason for the numerous 
works in this research line is that 
there are numerous interesting structures of graphs, and 
enumerating 
graphs satisfying certain structural constraints is a rich source of important 
research
questions, including enumerating trees, graphs with connectivity properties, 
graphs with fixed degree sequences, graphs with certain chromatic numbers, and so 
on. 

Therefore, a natural question is whether there are also interesting structures 
and properties of alternating matrix spaces. By far, alternating matrix spaces 
have received much less 
attention compared to graphs. But there is a classical connection between graphs 
and 
alternating matrix spaces,
going back to Tutte \cite{Tut47} and Lov\'asz \cite{Lov79}. This connection leads 
to interesting structures for 
alternating matrix spaces, which in turn forms the basis of our investigation.

\subsection{Some structures of alternating matrix spaces}\label{subsec:struct}

So let us review the classical construction of alternating matrix spaces from 
graphs \cite{Tut47,Lov79}.
Since we focus on labelled graphs, without loss of 
generality let us consider graphs with vertex sets being $[n]=\{1, \dots, n\}$. 
Let $G=([n], E)$ be an undirected simple graph, so $E\subseteq\binom{[n]}{2}$. 
For 
$\{i,j\}\in \binom{[n]}{2}$, $i<j$, an elementary alternating matrix $A_{i,j}\in 
\M(n, \F)$ is the $n\times n$ matrix with the $(i,j)$th entry being $1$, the 
$(j,i)$th entry being $-1$, and the rest entries being $0$. We then define 
\begin{equation}\label{eq:bA_G}
\cA_G=\langle A_{i,j} : \{i,j\}\in E \rangle\leq\Lambda(n, 
\F),
\end{equation}
where $\langle \cdot\rangle$ denotes linear span.
As $G$ has a perfect matching if and only if $\cA_G$ 
contains 
a full-rank matrix, Tutte used $\cA_G$ to 
characterise graphs 
without perfect 
matchings \cite{Tut47}, and Lov\'asz used $\cA_G$ to obtain a simple efficient 
randomised 
algorithm for the perfect matching problem \cite{Lov79}. 

%As we have seen, the construction of alternating matrix spaces from graphs in 
%Equation~\ref{eq:bA_G} leads to a 
%correspondence between perfect matchings, a graph-theoretic structure, and 
%full-rank matrices, a linear-algebraic structure. 
Inspired by this classical 
example, 
several correspondences between 
graph-theoretic structures, and structures for alternating matrix spaces, have 
been 
 discovered recently, including:
\begin{enumerate}
\item Independent sets vs isotropic spaces, and vertex colourings vs isotropic 
decompositions \cite{BCG+19};
%\item  \cite{BCG+19};
\item 
Connectivity vs orthogonal indecomposability. As a consequence, 
correspondences of vertex and edge connectivities for alternating matrix spaces 
are also presented \cite{LQ19};
\item Isomorphism notions for graphs and alternating matrix spaces \cite{HQ20a}.
\end{enumerate}
For now, let us review the one between 
independent sets and isotropic spaces in \cite{BCG+19}. Some other 
correspondences will be introduced later, when we come to the relevant enumeration 
problems. 
\begin{definition}\label{def:iso}
Let $\cA\leq\Lambda(n, \F)$ be an alternating matrix spaces. A subspace $U\leq 
\F^n$ is called a \emph{totally-isotropic space} of $\cA$, if for any $u, u'\in 
U$, and 
any 
$A\in \cA$, $u^tAu'=0$. 
\end{definition}
The correspondence between independent sets and totally-isotropic spaces is 
supported by 
the following.
Recall that $\alpha(G)$ denotes \emph{independence number} of a graph $G$.
%, 
%$\alpha(G)$, is the maximum size over all independent sets of $G$. 
Similarly, define the  
\emph{totally-isotropic number} of $\cA$, $\alpha(\cA)$, to be the maximum 
dimension over 
all totally-isotropic spaces of $\cA$. Letting $\cA_G$ be constructed from a graph 
$G$ as in Equation~\ref{eq:bA_G}, it is shown in 
\cite{BCG+19} that $\alpha(G)=\alpha(\cA_G)$. Based on this correspondence, 
several classical questions for independent sets are found to have natural 
correspondences in the alternating matrix space setting, with applications to 
group theory and quantum information theory \cite{BCG+19}.

\subsection{Overview of our results}\label{subsec:overview}

We now give an overview of our results. 
It is well-known that when setting $q=1$, the 
Gaussian binomial coefficient $\qbinom{n}{d}{q}$ becomes the normal binomial 
coefficient $\binom{n}{d}$. Almost all our results also have such a nice 
feature: when setting $q=1$, they become the corresponding graph enumeration 
formulas. 

For $c$-colorable graphs, there is a classical formula by Read \cite{Rea60} as 
reviewed in Equation~\ref{eq:read}. As a graph could be $k$-coloured in several 
ways, that formula essentially enumerates $(G, U)$, where $G=([n], E)$ is a graph, 
and $U$ is a partition of $[n]$ such that each subset in $U$ is an independent 
set of $G$. In 
Section~\ref{sec:iso}, we present a $q$-analogue of Read's formula 
for alternating matrix spaces with totally-isotropic $c$-decompositions (cf. 
Definition~\ref{def:isod}) in
Equation~\ref{eq:read_q}. 

For connected graphs, there is a classical formula as 
reviewed in Equation~\ref{eq:con}. Enumerating connected graphs is one of the 
first few questions studied in graph enumeration, tracing back to at least  
Riddell's thesis \cite{Rid51}, though it seems to the author that the exact 
formula as in Equation~\ref{eq:con} first appeared in Gilbert's article 
\cite{Gil56}. 

For alternating matrix spaces, there are two indecomposability notions, which, 
when restricted to alternating matrix spaces of 
the form $\cA_G$ as in Equation~\ref{eq:bA_G}, both coincide with the connectivity 
of $G$ \cite{LQ19}. 

The first indecomposability notion is with respect to the so-called direct 
decompositions (cf. Definition~\ref{def:ortho}). This indecomposability has a nice 
property that makes it quite close to connectivity. That is, for non-degenerate 
alternating matrix spaces (cf. Section~\ref{sec:prel}), there is a unique complete 
direct 
decomposition \cite{Wil12}. In Section~\ref{sec:free_ortho}, 
we study enumerating directly-indecomposable alternating matrix spaces. Despite 
some minor subtleties caused by the non-degeneracy condition, we obtain a 
$q$-analogue of Gilbert's formula in 
Equation~\ref{eq:di}.

The second indecomposability notion is with respect to the so-called orthogonal 
decompositions (cf. Definition~\ref{def:ortho}). This indecomposability is more 
flexible, meaning an alternating matrix space could have several complete
orthogonal decompositions in several ways even under the automorphism group action 
\cite{Wil09a}. So we give a formula for 
alternating matrix spaces with orthogonal $c$-decompositions in 
Equation~\ref{eq:ortho} in the spirit of Equation~\ref{eq:read_q}.

A basic tool in labelled graph enumeration is the exponential generating function. 
A very useful lemma based on exponential generating functions is the labelled 
counting lemma \cite[pp. 8]{HP73}. We demonstrate a $q$-analogue of this lemma for 
alternating matrix spaces, which we call the coordinate-explicit counting lemma, 
in Lemma~\ref{lem:ce}. As an application of this lemma, we derive a $q$-analogue 
of Riddell's formula relating the exponential generating function of 
graphs with that of connected graphs \cite{Rid51} in Equation~\ref{eq:rid_q}. 
Here, we heavily rely on the Eulerian generating 
functions developed by Srinivasan \cite{Sri06}.

%In the following, we shall study several enumerating questions for alternating 
%matrix spaces with explicit coordinates. As the reader will see, the results are 
%not difficult, though certain subtleties are to be addressed. Rather, the purpose 
%of these results is to support Postulate~\ref{main}.
%
%As the reader will see, the results in 
%Section~\ref{sec:iso} and~\ref{sec:free_ortho} can be viewed as nice linear 
%algebraic analogues, or $q$-analogues, of certain classical results in 
%enumerating 
%labelled graphs. In Section~\ref{sec:lem} we derive a linear algebraic analogue 
%of 
%the useful labelled counting lemma.
%
%Then in 
%Section~\ref{sec:ortho}, we will study a structure for alternating matrix spaces 
%that, in some sense, do not correspond to graph-theoretic structures. Still, 
%because of the explicit coordinates, we will manage to give a solution. 

\subsection{Remarks for future research}

We believe that the results in Section~\ref{subsec:overview} provide some initial 
support to our Postulate~\ref{main}. Of course, these results are mostly 
straightforward, and the corresponding graph enumeration results were all known in 
the 1950's. It is our hope that more $q$-analogues of formulas from graph 
enumeration will be found in the near future. 

For graph enumeration, the usefulness of the formulas is usually evidenced by the 
ability to calculate the exact numbers of graphs of order $n$ satisfying certain 
properties (as seen in e.g. \cite{HP73}). On the 
contraty, we do not expect our formulas could be used so, as the Galois numbers 
(see Section~\ref{sec:prel}) are already difficult to be evaluated exactly 
\cite{GR69}. Our main goal, at present, is to derive $q$-analogues of some exact 
counting formulas in graph enumeration, to support Postulate~\ref{main} and to get 
some aesthetic pleasure. 

In the future, we expect that, instead of exact recursive
enumeration formulas, asymptotic enumerations could be more useful for e.g. 
algorithms and 
probabilistic analysis, as also in the case of 
graph enumeration. Previous works of the author with collaborators 
\cite{LQ17} (improved later in \cite{BLQW}) support this possibility. Let $c$ be a 
constant. In 
\cite{LQ17}, it 
is shown that the automorphism group of a random $cn$-dimensional alternating 
matrix space $\cA$ in $\Lambda(n, q)$ is of order $q^{O(n)}$. Here, a random 
$\cA$ means drawing uniformly random from the 
$\qbinom{\binom{n}{2}}{m}{q}$-many $m$-dimensional subspaces of $\Lambda(n, q)$. 
This is a $q$-analogue of the well-known Erd\H{o}s-R\'enyi model \cite{ER1959}. 
Therefore the result in \cite{LQ17} can be viewed as the $q$-analogue of the 
celebrated result that when $m$ is asymptotically within $\frac{1}{2}n\log n$ and 
$\binom{n}{2}-\frac{1}{2}n\log n$, the automorphism group of a random graph with 
$n$ vertices and $m$ edges is 
trivial (i.e. of order $1$) \cite{ER63,Wri71}.

We didn't touch the topic of enumerating coordinate-free alternating matrix 
spaces, which corresponds to enumerating unlabelled graphs. Enumerating 
coordinate-free alternating matrix spaces can be understood as enumerating the 
orbits of $\Lambda(n, q)$ under the natural action of $\GL(n, q)$ (see 
Section~\ref{sec:prel}). We leave this topic to a future work. 

%Let us also indicate the differences between graph enumeration and alternating 
%matrix space numeration. 
%
%On the one hand, there are so many special and interesting graph families, like 
%trees, 
%tournaments, cubic graphs, and so on. Enumerating such graphs have been a rich 
%source of problems in graph enumeration. For alternating matrix spaces, despite 
%some recent discoveries mentioned in Section~\ref{subsec:struct}, special and 
%interesting families of alternating matrix spaces are not so abundant. 
%
%On the other hand, because of the ``correlation'' in the linear algebraic world, 
%new phenomena do arise in the alternating matrix space setting. In 
%Section~\ref{subsec:overview}, we have seen that connectivity can correspond to 
%two indecomposability notions, and the orthogonal decomposability notion does 
%exhibit new behaviours. 
%
%To summarise, the discussions in the above suggest that while an analogy can be 
%drawn between graph enumeration and alternating matrix space enumeration, we also 
%expect remarkable differences between these two research lines. 

\section{Preliminaries}\label{sec:prel}

%We collect some notation and useful results in this section. 

\paragraph{Some notions for alternating matrix spaces.}
For $\cA\leq\Lambda(n, \F)$ and $T\in\GL(n, \F)$, $T$ acts on $\cA$ naturally by 
sending $\cA$ to $T^t\cA T=\{T^tAT : A\in\cA\}$. Then the \emph{automorphism 
group} of $\cA$, 
$\aut(\cA)=\{T\in\GL(n, \F) : \cA=T^t\cA T\}$. We say that $\cA, \cB\leq\Lambda(n, 
\F)$ are isomorphic, if there exists $T\in\GL(n, \F)$, such that $\cA=T^t\cB T$.

For $\cA\leq\Lambda(n, \F)$, the \emph{radical} of $\cA$ is $\rad(\cA):=\{v\in\F^n 
: 
\forall A\in 
\cA, Av=0\}$. We say that $\cA$ is \emph{degenerate}, if $\rad(\cA)\neq 0$. 

Given a
dimension-$d$ $U\leq
\F^n$, let $T$ be a matrix of size $n\times d$ whose columns span $U$. The 
\emph{restriction} of $\cA$ on
$U$
via $T$ is $\cA|_{U, T}:=\{T^tAT : A\in \cA\}\leq \Lambda(d, \F)$. For a different 
$T'$ whose 
columns also span $U$, $\cA|_{U, T'}$ is isomorphic to $\cA|_{U, T}$. So we can 
write $\cA|_U$ to indicate a restriction of $\cA$ to $U$ via some such $T$. 

%if there exists a 
%non-zero $u\in \F^n$, such that for any $A\in \cA$, $v\in \F^n$, $u^tAv=0$.

\paragraph{Some basic $q$-calculus.} We present some notation and basic facts from 
the 
so-called $q$-calculus \cite{KC02}. For $n\in\N$ and a prime power $q$, let 
$[n]_q:=\frac{q^n-1}{q-1}=q^{n-1}+\dots+1$. 
The $q$-factorial $[n]_q!=[n]_q\cdot [n-1]_q\cdot \ldots\cdot [1]_q$.
The Gaussian binomial coefficient $\qbinom{n}{d}{q}$ can then be written as 
$\frac{[n]_q!}{[d]_q![n-d]_q!}$, which counts the number of dimension-$d$ 
subspaces of $\F_q^n$. 
%Expanding it, we see that 
%$$
%\frac{[n]_q!}{[d]_q![n-d]_q!}=
%\frac{q^{n-1}+\dots+1}{q^{d-1}+\dots+1}\cdot 
%\frac{q^{n-2}+\dots+1}{q^{d-2}+\dots+1}\cdot\ldots\cdot \frac{q^{n-d}+\dots+1}{1},
%$$
%so setting $q=1$ in the above, we have $\qbinom{n}{d}{1}=\binom{n}{d}$.

Following Goldman and Rota \cite{GR69} (see also 
\cite{KC02}), define the $n$th Galois number over $q$  
$\GalN{q}{n}$ to be $\sum_{d=0}^n\qbinom{n}{d}{q}$, which is
the total 
number of subspaces of $\F_q^n$. 
%By the discussion in the last paragraph, we have 
%$\GalN{1}{n}=2^n$. 

Let $\GN_n=2^{\binom{n}{2}}$ be the number of labelled graphs of order $n$. Let 
$\GN_{n, q}=\GalN{q}{\binom{n}{2}}$ be the number of alternating matrix spaces in 
$\Lambda(n, q)$. 

%\paragraph{Some correspondences between combinatorics and linear algebra.} Let 
%$S$ 
%be a finite set. We write a partition of $S$ as $S=S_1\uplus \dots\uplus S_c$. A 
%tuple of subsets of $S$, $(S_1, \dots, S_c)$, is an ordered partition, if 
%$S=S_1\uplus \dots \uplus S_k$. The number of ordered partitions $(S_1, \dots, 
%S_c)$ of  is counted by the 
%multimonomial coefficient $\binom{n}{n_1, \dots, 
%n_c}=\frac{n!}{n_1!\cdot\ldots\cdot n_c!}$.
%
%Let $V$ be a finite-dimensional vector space. We write a direct sum decomposition 
%of $V$ as $V=V_1\oplus\dots\oplus V_k$. A tuple of subspaces of $V$, $(V_1, 
%\dots, 
%V_k)$, is an ordered direct sum decomposition of $V$, if $V=V_1\oplus \dots\oplus 
%V_k$. The number of ordered direct sum dec

\paragraph{Eulerian generating functions.} We recall the notions of Eulerian 
generating functions, which are $q$-analogues of exponential generating functions. 
There are actually two versions of Eulerian generating functions. 
%Let $f, 
%g:\Z^+\to 
%\N$ be two functions. 

Let $U$ be a vector space, and $V\leq U$ be a subspace. We can consider the 
quotient space $U/V$, or complement subspaces of $V$ in $U$. 
%There is a unique 
%quotient 
%space $U/V$ of $U$ with respect to $V$. However, 
Note that
complement subspaces of $V$ in 
$U$, i.e. those $W\leq U$, $W\cap V=0$, and $\langle W, V\rangle=U$, are not 
unique. If $V$ is a dimension-$d$ subspace of $U=\F_q^n$, then there are 
$q^{d(n-d)}$-many complement subspaces of $V$ in $U$. 

The distinction between quotient spaces and complement spaces leads to two types 
of Eulerian generating functions. The Eulerian generating function for quotient 
spaces was first studied by Goldman and Rota in \cite{GR70}. To the best of our 
knowledge, the Eulerian 
generating function for complement spaces was first used explicitly by Srinivasan 
in \cite{Sri06}. We won't go into a detailed comparison between these two notions, 
but only indicate a key here: the coefficient of the Gaussian convolution for the 
Eulerian generating function in \cite{GR70} is $\qbinom{n}{d}{q}$, while that in 
\cite{Sri06} is $\qbinom{n}{d}{q}\cdot q^{d(n-d)}$, where the extra $q^{d(n-d)}$ 
is to take care of multiple complement subspaces. 

In this article we shall use the Eulerian generating function in \cite{Sri06}.

\section{From vertex $c$-colourings to isotropic $c$-decompositions}\label{sec:iso}

Totally-isotropic spaces for an alternating matrix space $\cA\leq\Lambda(n, \F)$ 
were 
defined in Definition~\ref{def:iso}, and results in \cite{BCG+19} indicate that 
totally-isotropic spaces 
can be studied as a linear algebraic analogue of independent sets. In 
\cite{BCG+19}, the following notion is also proposed.
\begin{definition}\label{def:isod}
Let $\cA\leq\Lambda(n, \F)$. A direct sum decomposition $\F^n=U_1\oplus 
\dots\oplus U_c$ is a \emph{totally-isotropic $c$-decomposition} of $\cA$, if each 
$U_i$ 
is a
totally-isotropic space of $\cA$. 
\end{definition}
In the following we shall only consider totally-isotropic $c$-decompositions which 
are 
non-trivial, i.e. none of $U_i$'s are the zero space. The totally-isotropic 
decomposition 
number of $\cA$, $\chi(\cA)$, is the smallest $c\in\N$ such that $\cA$ admits a 
totally-isotropic $c$-decomposition. Let $G$ be a graph and let $\cA_G$ be 
constructed 
from $G$ as in Equation~\ref{eq:bA_G}. It is shown in \cite{BCG+19} that 
$\chi(G)=\chi(\cA_G)$. Therefore, totally-isotropic decompositions can be viewed 
as a 
linear algebraic analogue of vertex colorings. 

Read gave a formula for the number of labelled $c$-coloured graphs \cite{Rea60}, 
generalising a result of Gilbert \cite{Gil56}. In this section, we shall review 
Read's formula first, and then present the result on enumerating alternating 
matrix spaces with isotropic $c$-decompositions. 

\paragraph{Review of enumerating $c$-colored labelled graphs.} Following Read 
\cite{Rea60}, to enumerate $c$-colored labelled graphs with $n$ 
vertices, we go through the following steps. It should be noted that a graph can 
be $c$-colored in several ways. So the following in fact enumerates pairs of the 
form $(G, U)$ where $G=([n], E)$ is a graph, and $U=(U_1, \dots, U_c)$ is 
an ordered partition of $[n]$, such that each $U_i$ is an independent set of $G$.
\begin{enumerate}
\item Enumerate ordered $c$-partition of $n$, i.e. $(n_1, \dots, n_c)$, 
$n_i\in\Z^+$, $\sum_{i\in[c]}n_i=n$.
\item Fix $(n_1, \dots, n_c)$, an ordered partition of $n$. Enumerate the number 
of ways to allocate $n_i$ vertices with color $i$. This is counted by 
the multinomial coefficient $\binom{n}{n_1, \dots, 
n_c}=\frac{n!}{n_1!\cdot\ldots\cdot n_c!}$.
\item Fix an allocation of $n_i$ vertices with color $i$. The number of graphs 
with this configuration is $2^{\binom{n}{2}-\sum_{i\in[c]}\binom{n_i}{2}}$, as 
only edges from one color class to another are possibly present.
\end{enumerate}
Summarising the above, Read's formula is 
\begin{equation}\label{eq:read}
\sum_{(n_1, \dots, n_c)}\binom{n}{n_1, \dots, 
n_c}\cdot 2^{\binom{n}{2}-\sum_{i\in[c]}\binom{n_i}{2}},
\end{equation}
where $(n_1, \dots, n_c)$, $n_i\in\Z^+$ goes over all ordered $c$-partition of $n$.
If the colors are not assume to have identity as in \cite{HP73}, then a 
multiplicative factor of $\frac{1}{c!}$ is required. 

\paragraph{Enumerating alternating matrix spaces with totally-isotropic 
$c$-decompositions.}
%
%Let $\cA\leq\Lambda(n, \F)$ be an alternating matrix space. We say that $V\leq 
%\F^n$ is a totally isotropic space of $\cA$, if for any $v, v'\in V$, $A\in \cA$, 
%it holds that $v^tAv'=0$. A direct sum decomposition $\F^n=V_1\oplus \dots \oplus 
%V_k$ is called a $k$-isotropic decomposition of $\cA$, if every $V_i$ is an 
%isotropic space of $\cA$. 
Following Read's recipe, we can enumerate alternating matrix spaces with 
totally-isotropic $c$-decompositions. 
\begin{enumerate}
\item Enumerate ordered $c$-partitions of $n$, i.e. $(n_1, \dots, n_c)$, 
$n_i\in\Z^+$, $\sum_{i\in[c]}n_i=n$. 
\item Fix $(n_1, \dots, n_c)$, an ordered partition of $n$. Enumerate the number 
of tuples of subspaces $(U_1, \dots, U_c)$, $U_i\leq \F_q^n$, such that 
$\dim(U_i)=n_i$, and $\F_q^n=U_1\oplus \dots \oplus U_c$. We shall refer to this 
tuple of subspaces an \emph{ordered direct sum decomposition} of $\F_q^n$. This 
number 
is, by 
\cite[Lemma 4]{Sri06},
$$
\frac{q^{\binom{n}{2}}\cdot [n]_q!}{(q^{\binom{n_1}{2}}\cdot [n_1]_q!)\dots 
(q^{\binom{n_c}{2}}\cdot [n_c]_q!)}.
$$
\item Fix an ordered direct sum decomposition $(U_1, \dots, U_c)$ of $\F_q^n$. 
Requiring $U_i$'s where $\dim(U_i)=n_i$ to be totally-isotropic spaces impose 
$\sum_{i\in[c]}\binom{n_i}{2}$ independent linear conditions on $\Lambda(n, q)$. 
So the number 
of alternating 
matrix spaces with $U_i$'s being isotropic spaces is 
$\GalN{q}{\binom{n}{2}-\sum_{i\in[c]}\binom{n_i}{2}}$.
\end{enumerate}
Summarising, we obtain a $q$-analogue of Read's formula (Equation~\ref{eq:read}) 
which counts the number of alternating matrix spaces with totally-isotropic 
$c$-decompositions:
\begin{multline}\label{eq:read_q}
  \sum_{(n_1, \dots, n_c)}\frac{q^{\binom{n}{2}}\cdot 
[n]_q!}{(q^{\binom{n_1}{2}}\cdot [n_1]_q!)\dots 
(q^{\binom{n_c}{2}}\cdot [n_c]_q!)}\cdot 
\GalN{q}{\binom{n}{2}-\sum_{i\in[c]}\binom{n_i}{2}} \\
=
\sum_{(n_1, \dots, n_c)}\frac{
[n]_q!}{[n_1]_q!\dots 
[n_c]_q!}\cdot 
q^{\binom{n}{2}-\sum_{i\in[c]}\binom{n_i}{2}}\cdot 
\GalN{q}{\binom{n}{2}-\sum_{i\in[c]}\binom{n_i}{2}},
\end{multline}
where $(n_1, \dots, n_c)$ goes over all ordered $c$-partitions of $n$.
If normal direct sum decompositions rather than 
ordered direct sum decompositions are counted, then a 
multiplicative factor of $\frac{1}{c!}$ is required. 
A nice feature of Equation~\ref{eq:read_q} is that when setting $q=1$ there, it 
becomes Equation~\ref{eq:read}.
%
%To do that, first consider ordered direct sum decompositions, i.e. $(V_1, \dots, 
%V_k)$ such that $\F^n=V_1\oplus \dots\oplus V_k$. Obviously, the ordered number 
%is 
%$k!$ times the unordered number. Second, consider ordered direct 
%sum decompositions with a fixed dimension sequence, i.e. those $(V_1, \dots, 
%V_k)$ 
%such that $\dim(V_i)=n_i$. Third, count the ordered direct sum decompositions 
%with 
%a fixed dimension sequence by \cite[Lemma 4]{Sri06}, which is
%$$
%\frac{q^{\binom{n}{2}}\cdot [n]_q!}{(q^{\binom{n_1}{2}}\cdot [n_1]_q!)\dots 
%(q^{\binom{n_k}{2}}\cdot [n_k]_q!)}.
%$$
%Fourth, for a fixed ordered direct sum decomposition with the degree sequence 
%$(n_1, \dots, n_k)$, we need to count subspaces of a linear space of dimension 
%$N=\binom{n}{2}-\sum_{i\in[k]}\binom{n_i}{2}$, which is $A_{N, q}$. To summarise 
%the above, the number is 
%$$
%\frac{1}{k!}\cdot \sum_{k\text{-ordered partition of 
%}n}\frac{q^{\binom{n}{2}}\cdot [n]_q!}{(q^{\binom{n_1}{2}}\cdot [n_1]_q!)\dots 
%(q^{\binom{n_k}{2}}\cdot [n_k]_q!)}\cdot A_{N, q}.
%$$

\section{From connectivity to direct and orthogonal
indecomposabilities}\label{sec:free_ortho}

Let us first define the structures of alternating matrix spaces to be studied in 
this section. We say that a direct sum decomposition of $\F_q^n=U_1\oplus \dots 
\oplus U_k$ is non-trivial, if none of $U_i$ is the zero space.
\begin{definition}\label{def:ortho}
For an alternating matrix space $\cA\leq 
\Lambda(n, \F)$, a non-trivial direct sum decomposition $\F^n=U_1\oplus 
\dots\oplus U_k$ is an 
\emph{orthogonal 
decomposition} of $\cA$, if for any $i\neq j$, $u_i\in U_i$, $u_j\in U_j$, and 
$A\in \cA$, $u_i^tAu_j=0$. 
An orthogonal decomposition is a \emph{direct decomposition}, if 
$\dim(\cA)=\sum_{i\in[k]}\dim(\cA|_{U_i})$.

We say that $\cA$ is \emph{orthogonally decomposable}, if it admits an 
orthogonal decomposition into $k\geq 2$ subspaces. It is \emph{directly 
decomposable}, if it 
admits a direct decomposition into $k\geq 2$ subspaces. 

An orthogonal (resp. direct) decomposition $U=U_1\oplus\dots\oplus U_k$ is 
\emph{complete}, if 
for any $U_i$, $\cA|_{U_i}$ is orthogonally (resp. directly) indecomposable. 
\end{definition}

Let $G$ be a graph, and let $\cA_G$ be constructed 
from $G$ as in Equation~\ref{eq:bA_G}. In 
\cite{LQ19}, it was shown that $G$ is 
disconnected if and only if $\cA_G$ is orthogonally decomposable. Actually, it is 
straightforward to 
see that orthogonally decomposable there can be strengthened to directly 
decomposable.
Therefore, both orthogonal indecomposability and direct indecomposability can be 
viewed as linear algebraic analogues of connectivity. 

%
%\begin{remark}\label{rem:ind}
%In \cite{LQ19}, we define $\Lambda(1, \F)$ to be orthogonally 
%decomposable. Therefore, more precisely speaking, orthogonal indecomposability 
%corresponds to 
%$1$-connectivity, which excludes the graph with a single 
%vertex \cite[pp. 12]{Die17}. Here, we further define $\Lambda(1, \F)$ to be 
%directly decomposable, and this will be useful in the following.
%\end{remark}

Interestingly, orthogonal indecomposability and direct indecomposability behave 
quite differently. 
%To see this, let us define the following.
%\begin{definition}\label{def:complete}
%For $\cA\leq\Lambda(n, \F)$, a direct sum decomposition $\F^n=U_1\oplus 
%\dots\oplus U_k$ is a complete direct decomposition, if for any $u_i\in U_i$, 
%$u_j\in U_j$ 
%$\cA|_{U_i}$ is 
%directly indecomposable for every $U_i$. Complete orthogonal decompositions can 
%be 
%defined similarly. 
%\end{definition}
On one hand, by \cite[Lemma 6.9 (iii)]{Wil12}, if $\cA$ is 
non-degenerate, then $\cA$ has a unique complete direct decomposition. On the 
other hand, by 
\cite[Theorem 1.1 (ii)]{Wil09a}, the number of complete orthogonal decompositions 
can be any positive integer even under the automorphism group of $\cA$. 
Therefore 
it is 
not surprising that the same strategy for enumerating 
alternating matrix spaces with totally-isotropic $c$-decompositions can be applied 
to 
enumerating alternating matrix spaces with orthogonal $c$-decompositions in a 
straightforward fashion. So we present the formula here without giving the 
details:
\begin{equation}\label{eq:ortho}
  \sum_{(n_1, \dots, n_c)}\frac{q^{\binom{n}{2}}\cdot 
[n]_q!}{(q^{\binom{n_1}{2}}\cdot [n_1]_q!)\dots 
(q^{\binom{n_c}{2}}\cdot [n_c]_q!)}\cdot 
\GalN{q}{\sum_{i\in[c]}\binom{n_i}{2}}.
\end{equation}

In the following, we focus on direct decompositions. We first review Gilbert's 
formula for enumerating connected 
graphs in Section~\ref{subsec:con}. We then derive a $q$-analogue of this formula 
in the setting of enumerating directly indecomposable alternating matrix spaces in 
Section~\ref{subsec:con_q}. 

%We finally indicate a formula for enumerating 
%alternating matrix spaces with orthogonal $c$-decompositions in 
%Section~\ref{subsec:ortho_c}.

\subsection{Review of enumerating connected graphs}\label{subsec:con}

Let us review the classical formula of Gilbert \cite{Gil56} and its proof 
following 
Harary and Palmer \cite[Eq. 1.2.]{HP73}, which uses the 
notion of rooted graphs. Recall 
that a 
\emph{rooted graph} is a graph $G=([n], E)$ with one of the vertices $v\in [n]$, 
called 
the root, distinguished from others. 

Recall that $\GN_n=2^{\binom{n}{2}}$ denotes the number of labelled graphs of 
order $n$. Let 
$\RGN_n$ be the 
number of rooted, labelled graphs of order $n$. Let $\CGN_n$ be the number of 
connected labelled graphs of order $n$.
On one hand, a graph of order $n$ gives rise to $n$ rooted 
graphs, so $\RGN_n=n\cdot \GN_n$. On the other hand, we count the number of 
rooted, labelled 
graphs depending on the size of the connected component containing the root, as 
follows.
\begin{enumerate}
\item Enumerate $k\in [n]$ as the size of the connected component containing the 
root. 
\item Fix $k\in[n]$. Enumerate $S\subseteq [n]$, $|S|=k$. In the following, $S$ 
will 
contain the root. 
\item Fix $S\subseteq[n]$, $|S|=k$. Enumerate connected, labelled graphs with the 
vertex set being $S$. Enumerate $v\in S$, where $v$ is chosen as the root. 
\item Enumerate labelled graphs with the vertex set being 
$[n]\setminus S$.
\end{enumerate}
The above recipe gives that 
$
\RGN_n=\sum_{k=1}^nk\cdot \binom{n}{k}\cdot \CGN_k\cdot \GN_{n-k}.
$
Therefore by $\RGN_n=n\cdot \GN_n$, we get 
\begin{equation}\label{eq:con}
\CGN_n=\GN_n-\frac{1}{n}\cdot \sum_{k=1}^{n-1}k\cdot \binom{n}{k}\cdot \CGN_k\cdot 
\GN_{n-k}.
\end{equation}

\subsection{Enumerating directly indecomposable alternating matrix 
spaces}\label{subsec:con_q}

Following the recipe in Section~\ref{subsec:con}, we can enumerate directly 
indecomposable alternating matrix spaces as follows. However, a little care will 
be needed to handle degenerate alternating matrix spaces. 

\paragraph{The issue with degenerate alternating matrix spaces.} Let us first 
formally cite a result of J. B. Wilson, which forms the basis for our 
enumeration formula. Recall the notion of complete direct decompositions in 
Definition~\ref{def:ortho}.
\begin{theorem}[{\cite[Lemma 6.9 (iii)]{Wil12}}]\label{thm:unique}
Let $\cA\leq\Lambda(n, q)$ be non-degenerate. Then there exists a unique 
complete direct decomposition for $\cA$.
\end{theorem}
If $\cA\leq \Lambda(n, q)$ is degenerate, then it can have several complete 
direct 
decompositions, which leads to over-counting if we follow the recipe in 
Section~\ref{subsec:con} directly.
%For example, suppose $\cA=\{ \begin{bmatrix}
%A' & 0 \\
%0 & 0 
%\end{bmatrix} : A'\in\Lambda(n-1, q)\}\leq\Lambda(n, q)$. Then any direct sum of 
%the form $\F^n=U\oplus \langle e_n\rangle$ where $U$ is a complement subspace of 
%$\langle e_n\rangle$ is a complete direct decomposition of $\cA$. As a result, it 
%can lead to over-counting for enumerative purposes. 
Let us present an example to explain how this over-counting occurs. In 
analogy with rooted graphs, we define rooted alternating 
matrix spaces. That is, a \emph{rooted alternating matrix space} is a pair $(\cA, 
v)$ 
where $\cA\leq \Lambda(n, \F)$, and $v\neq 0 \in 
\F^n$ is called the root.
\begin{example}
Let us follow the recipe in Section~\ref{subsec:con} for $\Lambda(5, q)$. In Step 
1, set  
$k=2$. In Step 2, it is natural to enumerate direct sum 
decompositions $\F_q^5=U_1\oplus U_2$, where $\dim(U_1)=2$ and 
$\dim(U_2)=3$. Let $e_i$ be the $i$th standard 
basis vector, and $C=\begin{bmatrix}
0 & 1 & 0 \\
-1 & 0 & 0 \\
0 & 0 & 0
\end{bmatrix}$.
\begin{enumerate}
\item Let $U_1=\langle e_1, e_2\rangle$, $U_3=\langle e_3, e_4, 
e_5\rangle$, and the root vector be $e_2$. For $\cA|_{U_1}$ there is only one 
choice, namely $\Lambda(2, q)$. 
Then suppose $\cA|_{U_2}=\langle C\rangle$. %, i.e. $\cA|_{U_2}$ is degenerate. 
%So $\cA$ is of dimension 
%$2$.
\item Now consider $V_1=\langle e_1+e_5, e_2\rangle$, $V_2=\langle e_3, e_4, 
e_5\rangle$, and the root vector be $e_2$. Let $\cB|_{V_1}=\Lambda(2, q)$, and let 
$\cB|_{V_2}=\langle C\rangle$. 
\end{enumerate}
Clearly, $\cA$ and $\cB$ are the same alternating matrix space 
${\small \Big\langle \begin{bmatrix}
0 & 1 & 0 & 0 & 0 \\
-1 & 0 & 0 & 0 & 0 \\
0 & 0 & 0 & 0 & 0 \\
0 & 0 & 0 & 0 & 0 \\
0 & 0 & 0 & 0 & 0 
\end{bmatrix}, \begin{bmatrix}
0 & 0 & 0 & 0 & 0 \\
0 & 0 & 0 & 0 & 0 \\
0 & 0 & 0 & 1 & 0 \\
0 & 0 & -1 & 0 & 0 \\
0 & 0 & 0 & 0 & 0 
\end{bmatrix}\Big\rangle},$
but it will be 
counted twice when enumerating the direct sum decompositions $U_1\oplus U_2$ and 
$V_1\oplus V_2$. 
\end{example}

The above examples leads us to ponder on Gilbert's formula in 
Equation~\ref{eq:con}. 
The formula's correctness is self-evident, based on the obvious fact that for any 
graph, there is a unique way to partition the vertex set into a disjoint union of 
connected components. In the alternating matrix space setting, since 
Theorem~\ref{thm:unique} ensures the uniqueness of complete direct decompositions 
only in the non-degenerate case, it will be natural to focus on non-degenerate 
alternating matrix spaces. 

\paragraph{The number of non-degenerate alternating matrix spaces.}  
%fix this  
%issue caused by degeneracy, let us focus on non-degenerate alternating matrix 
%spaces. 
The number 
of non-degenerate alternating matrix spaces is in complete analogy with the number 
of graphs without isolated vertices. 
Let us first recall a basic formula for graphs without isolated vertices. Let 
$\NDG_n$ be the number of labelled graphs on $n$ vertices without isolated 
vertices. Set $\NDG_0=1$. We can relate $\NDG_n$ with $\GN_n$ as follows. 
First, fix a 
size-$k$ subset $S$
of $[n]$ as isolated vertices. Second, put a graph with $n-k$ vertices without 
isolated vertices on $[n]\setminus S$. This gives that
$
\GN_n=\sum_{i=0}^n\binom{n}{k}\NDG_{n-k}=\sum_{i=0}^n\binom{n}{k}\NDG_{k}.
$

Analogously, let $\NDS_{n,q}$ be the number of non-degenerate alternating matrix 
spaces in $\Lambda(n, q)$. Set $\NDS_{0, q}=1$. We can related $\NDS_{n, q}$ with 
$\GN_{n, q}$ as 
follows. First fix a dimension-$k$ subspace $U$ of $\F_q^n$ as the radical 
(defined in Section~\ref{sec:prel}), and 
choose any complement subspace $V$ of $U$ in $\F_q^n$. Second, put a 
non-degenerate 
alternating 
matrix space in $\Lambda(n-k, q)$ with the \emph{support space} being $V$. That 
is, first fix 
a linear isomorphism $T:\F_q^{n-k}\to V$, represented by a $(n-k)\times n$ matrix, 
and then send $\cA\leq \Lambda(n-k, q)$ to 
$T^t\cA T\leq\Lambda(n, q)$. This gives that
\begin{equation}\label{eq:NDS}
\GN_{n, q}=\sum_{i=0}^n\qbinom{n}{k}{q}\NDS_{n-k, 
q}=\sum_{i=0}^n\qbinom{n}{k}{q}\NDS_{k, q}.
\end{equation}
Because of Equation~\ref{eq:NDS}, we shall assume that $\NDS_{n, q}$ is known from 
$\GN_{n, q}$.
%, and 
%focus on non-degenerate alternating matrix spaces in the rest of this subsection. 

\paragraph{A $q$-analogue of Gilbert's formula.} Let $\DIS_{n, q}$ be the number 
of direct-indecomposable, non-degenerate alternating matrix 
spaces 
in $\Lambda(n, q)$. Let $\RSN_{n, q}$ be the number of rooted, non-degenerate 
alternating matrix 
spaces in $\Lambda(n, q)$. Recall that $\NDS_{n, q}$
denotes the number of non-degenerate
alternating matrix spaces in $\Lambda(n, q)$.
On one hand, an alternating matrix space in $\Lambda(n, q)$ yields $(q^n-1)$ 
rooted alternating matrix spaces, so $\RSN_{n,q}=(q^n-1)\cdot \NDS_{n,q}$. On the 
other 
hand, we count the number of rooted, non-degenerate alternating matrix spaces 
depending on the 
dimension of the subspace containing the root. 
\begin{enumerate}
\item Enumerate $k\in\{2, \dots, n\}$ as the dimension of the subspace containing 
the 
root.
\item Fix $k\in[n]$. Enumerate $U\leq \F^n$, $\dim(U)=k$, and enumerate complement 
subspaces of $U$. In the following, $U$ 
will contain the root.
\item Fix $U\leq \F^n$, $\dim(U)=k$. Enumerate directly indecomposable 
non-degenerate alternating 
matrix spaces with the support space being $U$. Enumerate non-zero vectors in 
$U$ as the root. 
\item Fix a complement subspace of $U$ in $\F^n$. Let it be $V$. Enumerate 
non-degenerate alternating matrix spaces with the support space being $V$. 
\end{enumerate}
The above recipe gives that
$
\RSN_{n, q}=\sum_{k=2}^n (q^k-1)\cdot \qbinom{n}{k}{q}\cdot q^{k(n-k)}\cdot 
\DIS_{k, q}\cdot 
\NDS_{n-k, q}.
$
Therefore using $\RSN_{n, q}=(q^n-1)\cdot \NDS_{n, q}$, we have
\begin{equation}\label{eq:di}
\DIS_{n, q}=\NDS_{n, q}-\frac{1}{[n]_q}\cdot\sum_{k=2}^{n-1}[k]_q\cdot 
\qbinom{n}{k}{q}\cdot q^{k(n-k)}\cdot \DIS_{k,q}\cdot \NDS_{n-k, q}.
\end{equation}

For graphs, it is also easy to derive a version of Gilbert's formula in 
Equation~\ref{eq:con} for graphs without isolated vertices. That formula would be 
the same as setting $q=1$ there in Equation~\ref{eq:di}.

\section{From labelled counting lemma to coordinate-explicit counting 
lemma}\label{sec:lem}

The exponential generating function and the labelled counting lemma are basic 
tools in graph enumeration. We first review them below. We then present a 
corresponding lemma for enumerating alternating matrix spaces, building on the 
work of Srinivasan \cite{Sri06}.

\paragraph{Review of the labelled counting lemma.}
Given a function $f:\Z^+\to\N$ and a variable $x$, the exponential generating 
function for $f$ in $x$ is 
$\exp(f, x)=\sum_{n\in \Z^+} f(n)\cdot \frac{x^n}{n!}.$
When the variable $x$ is understood from the context, we may simply write $\exp(f, 
x)$ as $\exp(f)$.
Given $f_i:\Z^+\to\N$, $i\in[c]$, suppose $\prod_{i\in[c]}\exp(f_i, 
x)=\sum_{n\in\N}f(n)\cdot 
\frac{x^n}{n!}=\exp(f, x)$, where $f:\Z^+\to\N$. Then 
$f(n)=\sum_{(n_1, \dots, n_c)} 
\binom{n}{n_1, \dots, n_c}\cdot \prod_{i\in[c]}f_i(n_i),$
where $(n_1, \dots, n_c)$, $n_i\in\Z^+$, goes over all ordered $c$-partition of 
$n$. 
The above, when interpreting in the context of graphs, leads to the following.
\begin{lemma}[{Labelled counting lemma, see e.g. \cite[pp. 8]{HP73}}]
Let $f_i, f:\Z^+\to\N$ be from above. Suppose $f_i:\Z^+\to\N$ 
is the exponential generating function for the number of labelled graphs 
satisfying property $P_i$. Then $f(n)$ counts 
the number of tuples of labelled graphs $(G_1, \dots, G_c)$, such that $G_i$ 
satisfies $P_i$, the sum of orders of $G_i$ is $n$, and the union of the vertex 
sets of $G_i$ is $[n]$.
\end{lemma}

An immediate application of the labelled counting lemma is the following relation 
discovered by Riddell \cite{Rid51} (cf. \cite[pp. 8]{HP73}). Recall that $\GN_n$ 
counts the number of labelled graphs of order $n$, and $\CGN_n$ counts the number 
of connected labelled graphs of order $n$. With a little manipulation, the 
labelled counting lemma gives that 
$\exp(\GN)=\sum_{c=1}^\infty \exp(\CGN)^c/c!$, which can be recorded conveniently 
as 
\begin{equation}\label{eq:rid}
1+\exp(\GN)=e^{\exp(\CGN)}.
\end{equation}

\paragraph{A coordinate-explicit counting lemma.}
Following Srinivasan \cite{Sri06}, we define the following. Given a function 
$f:\Z^+\to\N$, a variable 
$x$, 
and a prime power 
$q$, the 
Eulerian generating function for $f$ in $x$ and $q$ is
$$
\exp_q(f, x)=\sum_{n\in\Z^+}f(n)\cdot\frac{x^n}{q^{\binom{n}{2}}\cdot [n]_q! }.
$$
We may omit $x$ and write $\exp_q(f, x)$ in the following. 

Given $f_i:\Z^+\to \N$, $i\in[c]$, suppose $\prod_{i\in[c]}\exp_q(f_i, 
x)=\sum_{n\in\Z^+}f(n)\cdot \frac{x^n}{q^{\binom{n}{2}}\cdot[n]_q! }=\exp_q(f, 
x)$. 
Then by \cite[Theorem 5]{Sri06},
$$
f(n)= \sum_{(n_1, \dots, n_c)}\frac{q^{\binom{n}{2}}\cdot 
[n]_q!}{(q^{\binom{n_1}{2}}\cdot [n_1]_q!)\dots 
(q^{\binom{n_c}{2}}\cdot [n_c]_q!)} \cdot \prod_{i\in[c]}f_i(n_i),
$$
where $(n_1, \dots, n_c)$, $n_i\in\Z^+$, goes over all ordered $c$-partition of 
$n$. 

The above, when interpreting in the context of alternating matrix spaces, leads to 
the following.
\begin{lemma}[Coordinate-explicit counting lemma]\label{lem:ce}
Let $f_i, f:\Z^+\to\N$ be from above. Suppose $f_i:\Z^+\to\N$ is
the Eulerian generating function for the number of alternating matrix spaces 
satisfying property $P_i$. Then 
$f(n)$ counts 
the number of tuples of coordinate-explicit alternating matrix spaces $(\cA_1, 
\dots, \cA_c)$, such that $\cA_i$ 
satisfies $P_i$, $\cA_i$ supported by $U_i\leq \F_q^n$ of dimension $n_i\in\Z^+$, 
$\sum_{i\in[c]}n_i=n$, and $\F^n=U_1\oplus \dots \oplus U_c$.
\end{lemma}
(For the notion of $\cA_i$ supported by $U_i$, see the discussion before 
Equation~\ref{eq:NDS}.)

Based on the above, we can derive a $q$-analogue of Riddell's formula. The 
following is essentially a consequence of \cite[Theorem 6]{Sri06}. Recall that 
$\NDS_{n, q}$ counts the number of non-degenerate alternating matrix spaces in 
$\Lambda(n, q)$, 
and $\DIS_{n, q}$ counts the number of directly indecomposable, non-degenerate 
alternating matrix 
spaces in $\Lambda(n, q)$. %Recall that by Remark~\ref{rem:ind}, $\DIS_{1, q}=0$. 
By Lemma~\ref{lem:ce}, $\exp_q(\DIS)^c$ is 
the Eulerian 
generating function 
of ordered $c$-tuples of directly indecomposable, non-degenerate alternating 
matrix spaces, whose support spaces form a direct sum decomposition of the 
underlying 
vector space. Then 
the coefficient of $x^n/(q^{\binom{n}{2}}\cdot [n]_q!)$ of $\exp_q(\DIS)^c/c!$ is 
the number of non-degenerate alternating matrix spaces in $\Lambda(n, q)$ whose 
complete direct decompositions have $c$ summands. It follows that 
$\exp_q(\NDS)=\sum_{c=1}^\infty \exp_q(\DIS)^c/c!$, 
which can be convenient recorded as 
\begin{equation}\label{eq:rid_q}
1+\exp_q(\NDS)=e^{\exp_q(\DIS)}.
\end{equation}
Let us also add that the correctness of the above reasoning follows from 
Theorem~\ref{thm:unique}. 

%In particular, suppose we remove the non-degenerate 
%condition, i.e. replace $\NDS_{n, q}$ with $\GN_{n, q}$, and replace $\DIS_{n, 
%q}$ 
%with $\widetilde\DIS_{n,q}$ which counts the number of directly indecomposable 
%spaces. Then $1+\exp_q(\GN)=e^{\exp_q(\widetilde\DIS)}$ does not hold, as one 
%degenerate alternating matrix spaces gives rise to several direct 
%decompositions, which cannot be accounted by multiplying $1/c!$, as that only 
%removes the redundancy caused by ordering the summands in the direct sum 
%decompositions.
%
%It is interesting to examine whether in the above reasoning, $\NDS_{n, q}$ could 
%be replaced with $\GN_{n, q}$, to related with $\tilde \DIS_{n, q}$ which counts 
%the number of directly indecomposable alternating matrix spaces. In particular, 

\bibliographystyle{alpha}
\bibliography{references}

\end{document}